\documentclass[12pt]{article}

\usepackage{graphicx}
\usepackage{float}

\begin{document}

\title{Distributions of differences of Riemann zeta zeros}

\maketitle
\begin{center}
\author{J. Takalo}

{Dpt. of Physical Sciences, Space Climate research unit, University of Oulu,
POB 3000, FIN-90014, Oulu, Finland\\
email:{jouni.j.takalo@oulu.fi}}

\end{center}

\begin{abstract}
We study distributions of differences of unscaled Riemann zeta zeros, $\gamma-\gamma^{'}$, at large. We show, that independently of the location of the zeros, their differences have similar statistical properties. The distributions of differences are skewed towards the nearest zeta zero, have local maximum of variance and local minimum of kurtosis at or near each zeta zero. Furthermore, we show that distributions can be fitted with Johnson probability density function, despite the value of skewness or kurtosis of the distribution.
\end{abstract}


\section{Introduction}

Montgomery (1973) conjectured that the paircorrelation between pairs of zeros of the Riemann zeta function (scaled to have unit average spacing) is \cite{Montgomery}

\begin{equation}
\label{eq:CUE}
	R_{2}\left(x\right) = 1\!-\!\left(\frac{sin(\pi\,x)}{\pi\,x}\right)^{2} .
\end{equation}

\begin{figure}
	\centering
	\includegraphics[width=0.8\textwidth]{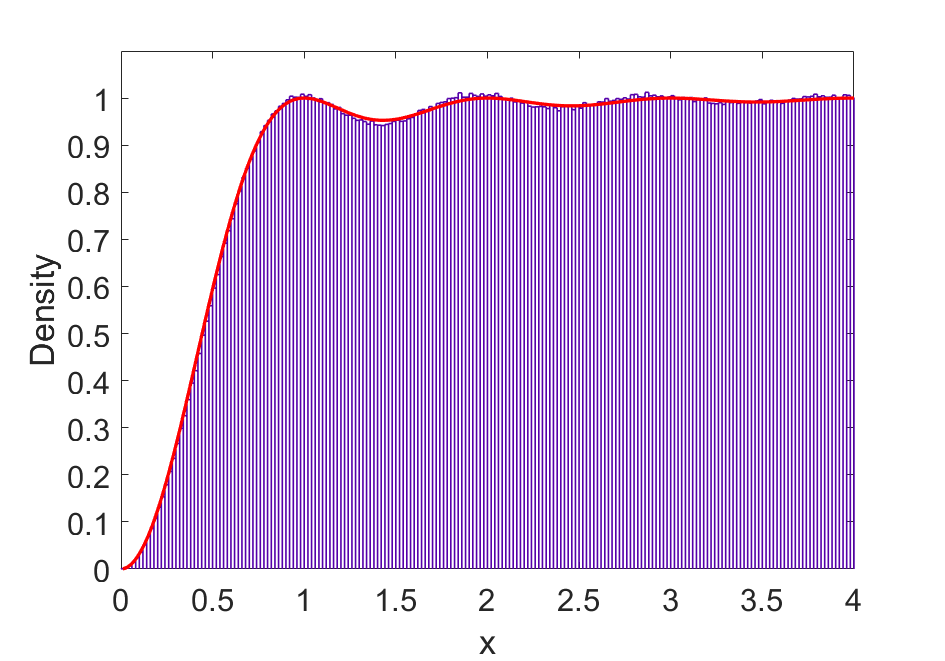}
		\caption{Distribution of five million consecutive and locally close scaled differences of zeta zeros starting from billionth zero. Red curve is the plot of function in Eq.1.}
		\label{fig:prediction}
\end{figure}

Figure \ref{fig:prediction} shows the density of 5 million differences starting from billionth zero and the prediction of Eq.1 as a red curve. Odlyzko (1987) already showed that the conjecture was supported by large-scale computer calculations of the zeros \cite{Odlyzko}. The aforementioned papers are restricted to analyse only consecutive or locally close differences of zeros. We study here distributions of differences of unscaled zeta zeros at large. Perez-Marco (2011) has shown that the statistics of very large zeros do find the location the first Riemann zeros. Our study differs from the others such that we analyse distributions of each difference $\gamma(n+i)-\gamma(i)$ separately for each n, and show that the properties of these separate distributions are the key to the large scale statistics of zeta zeros. Furthermore, we use Johnson probability density function (PCF) to show its flexibility for fitting zeta zero difference distributions.

\section{Data and methods}

\subsection{Nontrivial zeros of Riemann zeta function}

The data, i.e., the imaginary parts of nontrivial Riemann zeta zeros were fetched from  (https://www.lmfdb.org/zeros/zeta/) up to 100 billionth zero. The 10000 zeros at height $10^{21}$ and 10000 zeros at height $10^{22}$ are from \newline (www.dtc.umn.edu/~odlyzko/zeta\_tables/index.html)

\subsection{Johnson distribution}

Johnson distribution for the variable x is defined as 
\begin{equation}
z = \lambda+\delta\,ln\left(f\left(u\right)\right),
\end{equation}
with 
\begin{equation}
u = \left(x-\xi\right)/\lambda,
\end{equation}
Here z is a standardized normal variable and $f\left(u\right)$ has three different forms
the lognormal distribution, $S_{L}$:
\begin{equation}
f\left(u\right)=u,
\end{equation}
the unbounded distribution, $S_{U}$:
\begin{equation}
f\left(u\right)=u+{\left(1+u^{2}\right)}^{1/2},
\end{equation}
and the bounded distribution, $S_{B}$:
\begin{equation}
f\left(u\right)=u/\left(1-u\right).
\end{equation}
The supports for the distributions are $S_{L}: \xi<x,\; S_{U}: - \infty<x< \infty\; $and $S_{B}: \xi<x<\xi+\lambda$ \cite{Johnson, Wheeler}. The reason for the transformation of the the non-normal variables to standardized normal variables was that normal distribution was the only well-defined distribution at those times. However, with these definitions, the probability distributions are for
$S_{L}$:
\begin{equation}
P\left(u\right)=\frac{\delta}{\sqrt{2\pi}}\,\times\,\frac{1}{u}\,\times\,exp\left\{-\frac{1}{2}\left[\gamma+\delta\ln\left(u\right)\right]^{2}\right\}.
\end{equation}
for $S_{U}$
\begin{equation}
P\left(u\right)=\frac{\delta}{\sqrt{2\pi}}\,\times\,\frac{1}{\sqrt{u^{2}+1}}\,\times\,exp\left\{-\frac{1}{2}\left[\gamma+\delta\ln\left(u+\sqrt{u^{2}+1}\right)\right]^{2}\right\}.
\end{equation}
and for $S_{B}$
\begin{equation}
P\left(u\right)=\frac{\delta}{\sqrt{2\pi}}\,\times\,\frac{1}{u/\left(1-u\right)}\,\times\,exp\left\{-\frac{1}{2}\left[\gamma+\delta\ln\left(\frac{u}{1-u}\right)\right]^{2}\right\}
\end{equation}

\section{Distributions of differences of zeta zeros}

In the next we study the differences, $(\delta)$, of unscaled zeta zeros. We use the following notation

\begin{equation}
\label{eq:delta}
\delta(n) = \gamma(n+i)-\gamma(i), n=1,2,3,...
\end{equation}

\begin{figure}
	\centering
	\includegraphics[width=0.85\textwidth]{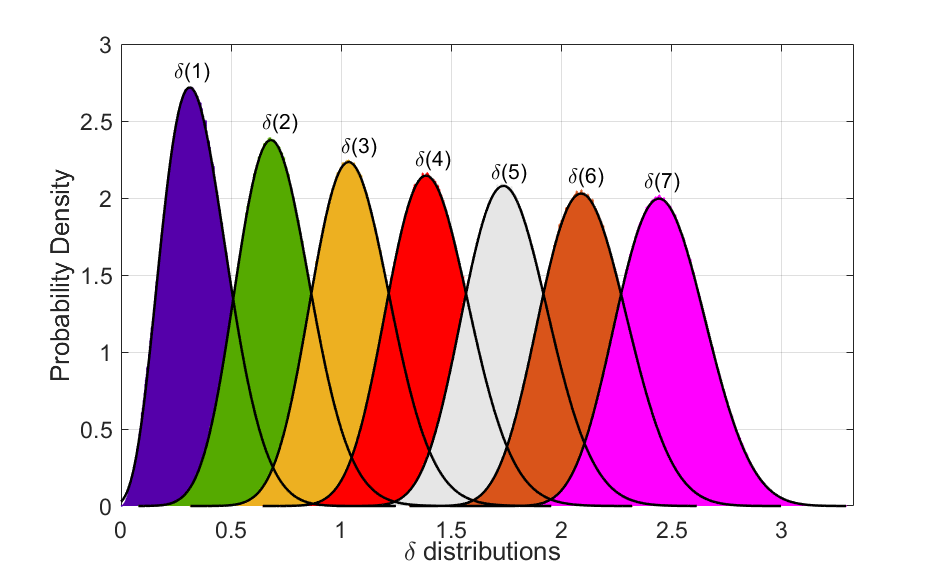}
		\caption{First seven $delta$-distributions calculated with 5 million unscaled zeta zeros starting from billionth zero.}
		\label{fig:First_zeta_1_billion_deltas}
\end{figure}

Here $i$ goes from j to j+5000000 in our analyses (if not otherwise mentioned), where j is the ordinal number of the starting zeta zero. Because the zeros are not stabilized yet at millionth zero we study, in the next, differences starting from 1, 10 and 100 billionth zeta zero, i.e., j=$10^{9}$, $10^{10}$ and $10^{11}$. 
Figure \ref{fig:First_zeta_1_billion_deltas} shows the first seven $delta$-distributions using five million differences starting from 1 billionth zeta zero. The distributions are fitted with the Johnson SB PDF. Notice, that the fits are quite good, except the fit of $delta(1)$ near zero. It is also evident that the height of the $delta$-distributions decrease and consequently the variance increases as a function of n. It is natural to ask how long this trend lasts.

\begin{figure}
	\centering
	\includegraphics[width=0.8\textwidth]{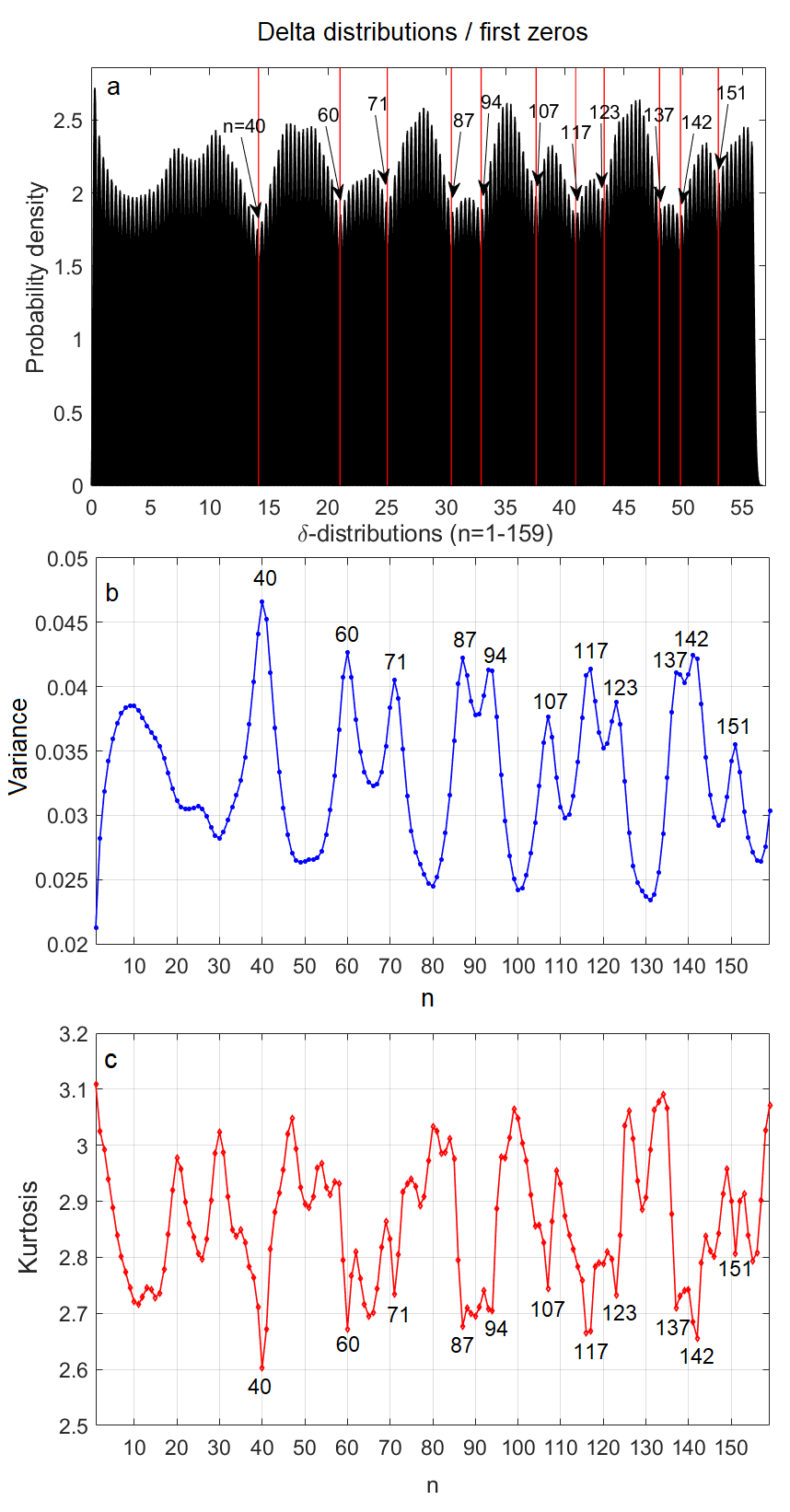}
		\caption{a) $Delta(n)$-distributions with n=1,2,...,159 calculated with 5 million unscaled zeta zeros starting from billionth zero.b) Variances of the distributions of figure a. c) Kurtoses of the distributions of figure a.}
		\label{fig:Zeta_1_billion_deltas_1_159}
\end{figure}

Figure \ref{fig:Zeta_1_billion_deltas_1_159}a shows the distributions of five million $deltas$ for n=1,2,3,..,159 from the same height as in Fig. \ref{fig:First_zeta_1_billion_deltas} together with 11 smallest zeta zero shown by red vertical line. Figure \ref{fig:Zeta_1_billion_deltas_1_159}b and \ref{fig:Zeta_1_billion_deltas_1_159}c show the corresponding variances and kurtoses of the 159 $delta$-distributions shown in Fig. \ref{fig:Zeta_1_billion_deltas_1_159}a, respectively. Notice that variances have local maxima and kurtoses local minima at distributions n=40, 60, 71, 87, 94, 107, 117, 123, 137, 142 and 151, i.e., at the sites, where distribution is nearest to each zeta zero. Note, that some zeros are between two consecutive distribution peaks, especially 95/96, and these distributions have almost the same variance/kurtosis between each other. The zeros 48.005 and 49.774 are so near each other that almost all distribution surrounding them, i.e., distributions 137-142 have almost similar variances and kurtoses. Note also that distribution 141 seems to have slightly higher variance than distribution 142, although the latter seems to be somewhat nearer to the nearby zero.  Kurtosis is, however, smaller for distribution 142 than distribution 141. The skewness is 0.0305 for 141 and -0.1042 for 142, and the higher absolute value of skewness decreases the variance of distribution 142. It is also notable that the skewness changes sign when passing the zeta zero.

\begin{figure}
	\centering
	\includegraphics[width=1.0\textwidth]{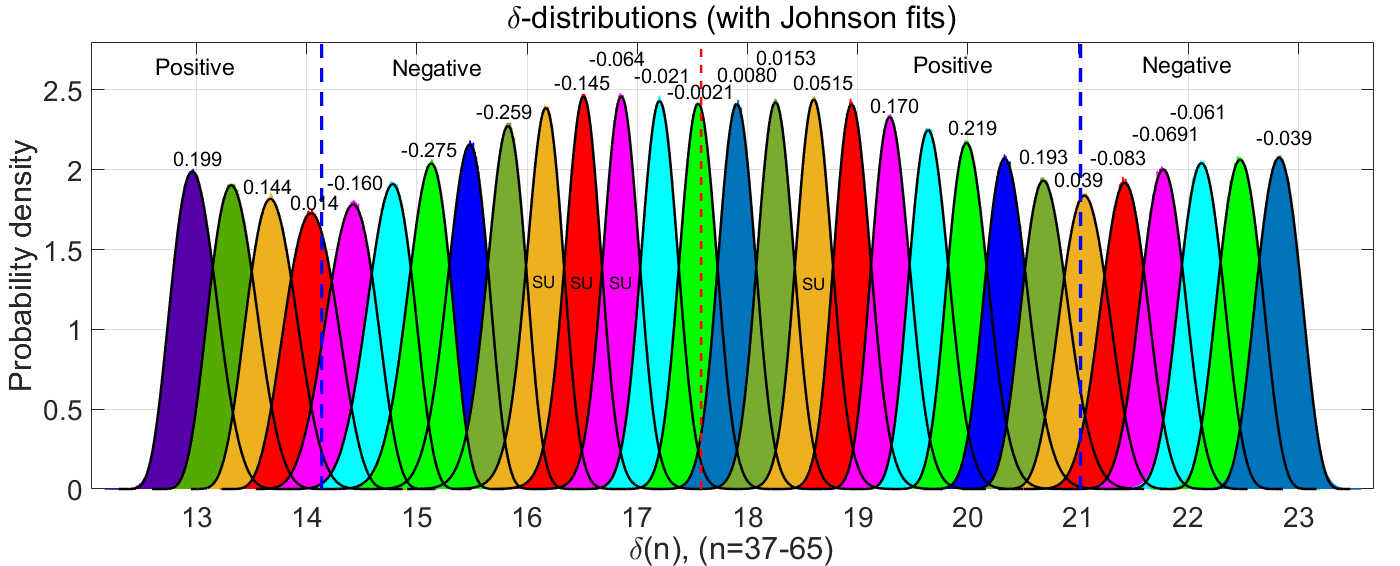}
		\caption{$Delta(n)$-distributions with n=37-65 fitted with Johnson probability density function. The decimal numbers are skewnesses of the separate distributions and the text tells their signs in each region.}
		\label{fig:More_detailed_1_billion_deltas_37_65}
\end{figure}

\begin{figure}
	\centering
	\includegraphics[width=0.85\textwidth]{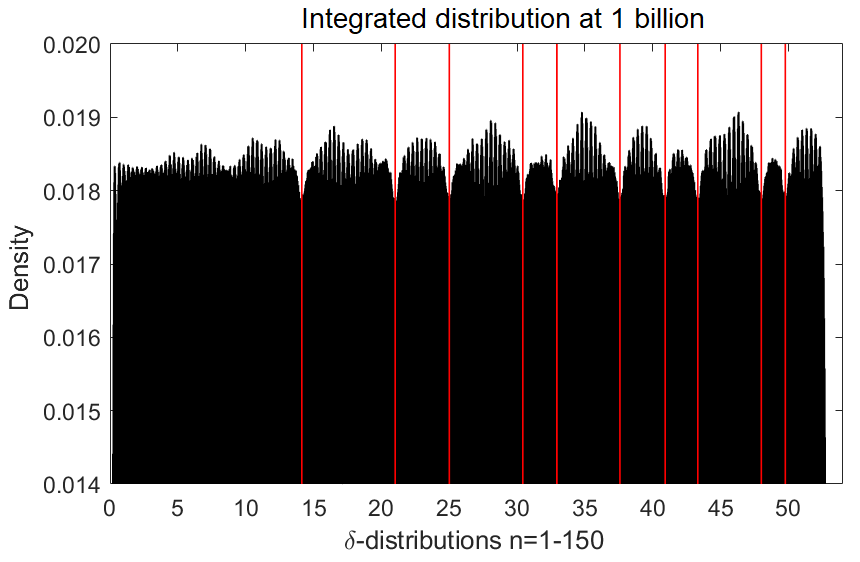}
		\caption{Integrated $delta(n)$-distributions of Fig. \ref{fig:Zeta_1_billion_deltas_1_159} for n=1,2,...,150.}
		\label{fig:Summed_1_billion_deltas_1_159}
\end{figure}

This is actually case around every zero. Figure \ref{fig:More_detailed_1_billion_deltas_37_65} shows more detailed pattern of the $delta$-distributions of n=37-65. The Johnson SB and SU (the latter marked in figure) fits are also shown in the figure. The blue dotted vertical lines show the two first zeta zeros, and red thinner dotted line the halfway between these zeros. The decimal number above each distribution is the skewness of corresponding distribution, and their sign in each region is also shown in the text above them. Note that the skewness changes sign at the zeros, and furthermore, at the halfway of the zeros such that $delta$-distribution is always skewed towards nearest zeta zero. It seems that the $deltas$ of the zeros are giving way for the zeros themselves, i.e., Riemann zeros do repel their $deltas$ \cite{Perez-Marco}. All the above-mentioned reasons lead to troughs around the zeros at the zeta zeros in integrated $delta$-distribution shown in figure \ref{fig:Summed_1_billion_deltas_1_159} \cite{Snaith, Perez-Marco}.

\begin{figure}
	\centering
	\includegraphics[width=0.85\textwidth]{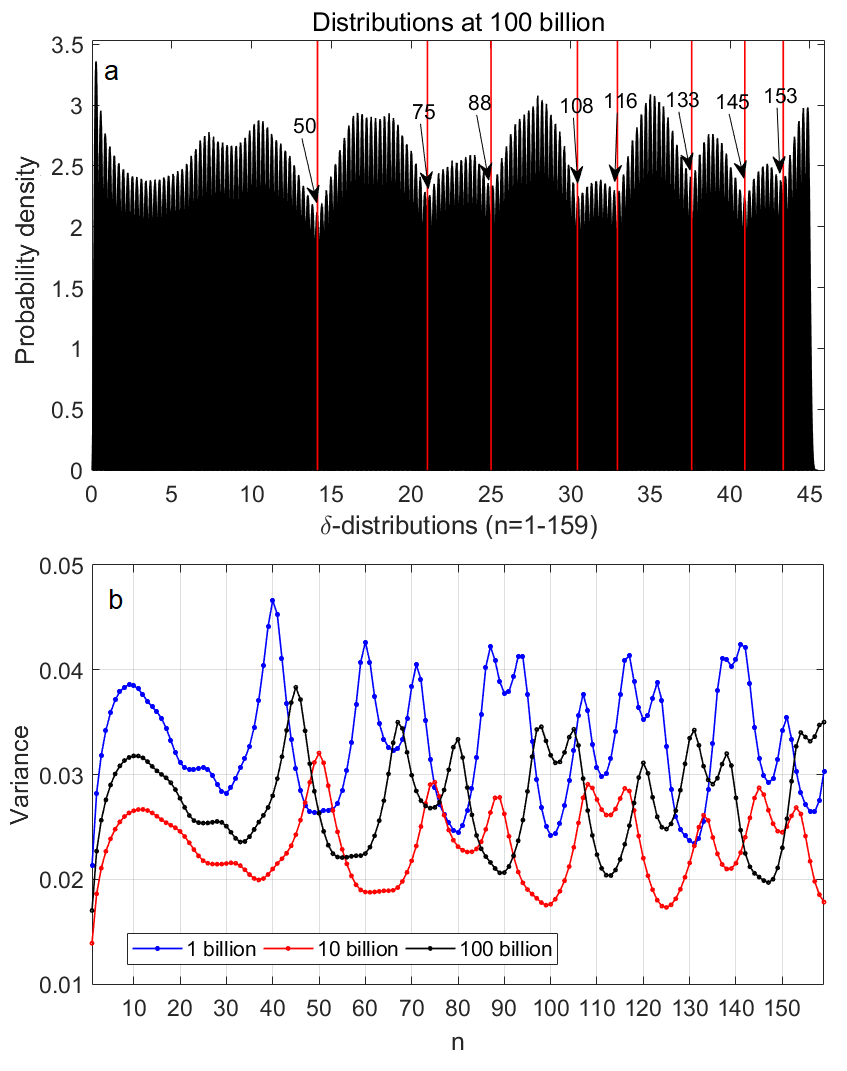}
		\caption{a) $Delta(n)$-distributions with n=1,2,...,159 calculated with 5 million unscaled zeta zeros starting from 100 billionth zero.b) Variances of the distributions of 5 million $deltas$ starting at 1 billionth (blue) 10 billionth (black) and 100 billionth (red) zero.}
		\label{fig:Zeta_100_billion_deltas_1_159}
\end{figure}

Figure \ref{fig:Zeta_100_billion_deltas_1_159}a shows the $delta$-distributions for 5 million differences with n=1-159 starting from 100 billionth zeta zero. The distributions behave similarly to that of Fig. \ref{fig:Zeta_1_billion_deltas_1_159}a for zeros at one billion. The only difference is that sites of zeros have moved forward such that they are now at $delta$-distributions n=50, 75, 88, 108, 116 133, 145 and 153. This is because variances are smaller and consequently the distributions at 100 billion height are narrower. Note, that the last three zeros of Fig. \ref{fig:Zeta_1_billion_deltas_1_159} are not seen anymore with the same number of distributions. Figure \ref{fig:Zeta_100_billion_deltas_1_159}b shows variances for $delta$-distributions at 1, 10 and 100 billionth zeta zeros. The pattern of the variances are similar, except that they are stretched for 10 billion and still more for 100 billion. The stretching ratios are about 1.125  between 1 billion and 10 billion and 1.25 between 1 billion and 100 billion.

\begin{figure}
	\centering
	\includegraphics[width=0.8\textwidth]{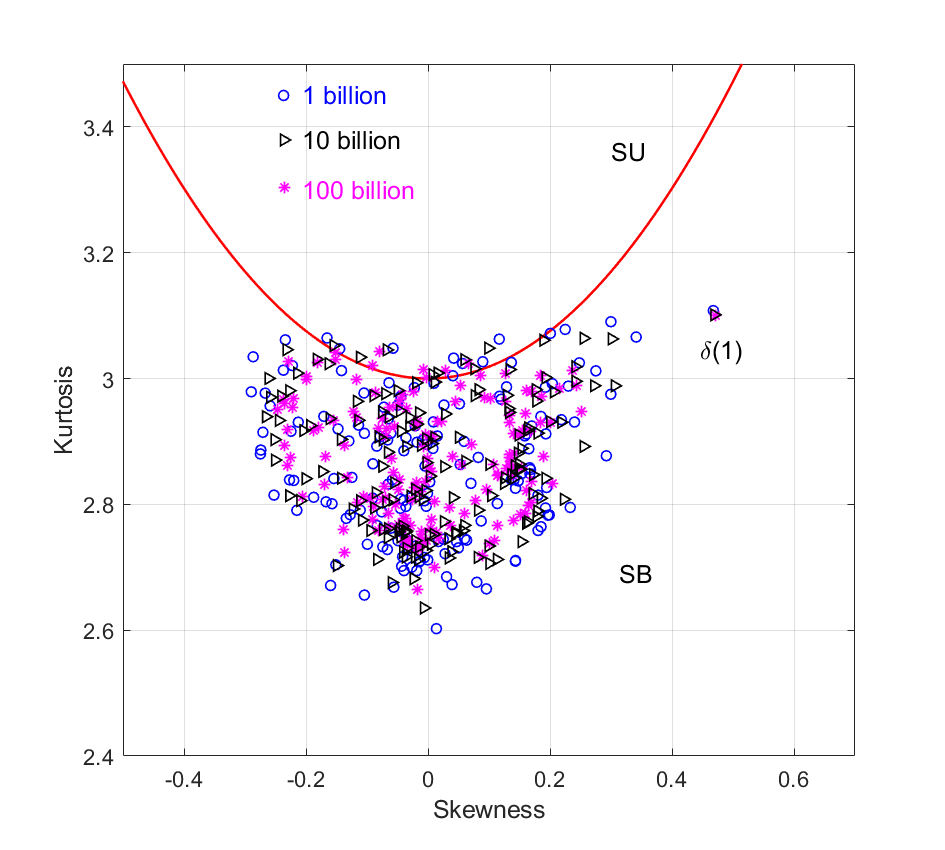}
		\caption{Skewness-kurtosis -plane with points of distributions starting 1 billionth (blue circle), 10 billionth (black triangle) and 100 billionth (magenta star) zero.}
		\label{fig:Skewness_kurtosis_coord}
\end{figure}

As shown earlier the $delta$-distributions can be fitted very well with the Johnson SB and SU PDFs. Figure \ref{fig:Skewness_kurtosis_coord} shows the distributions in the skewness-kurtosis -plane \cite{Cugerone}. The border between Johnson SU and SB distributions is marked with red curve in the plane. It is notable that the $delta(1)$s are in the same place in the plane. Otherwise, the points seem to be somewhat more compactly located when the level of zeta zeros increases.

\begin{figure}
	\centering
	\includegraphics[width=0.8\textwidth]{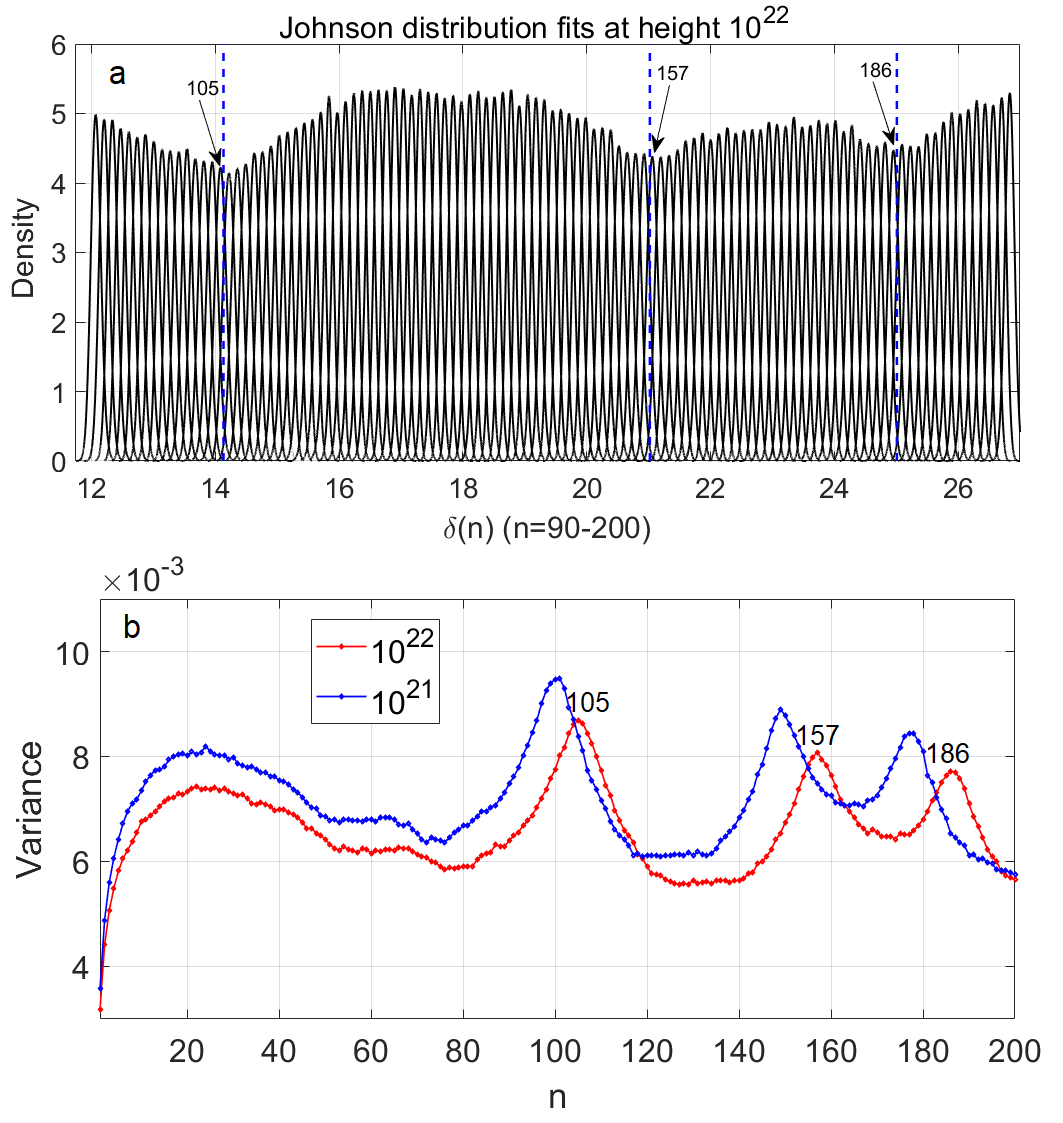}
		\caption{a) Johnson probability density fits for $delta(n)$ with n=1,2,3,...,200 at height $10^{22}$. b) Variances of the distributions at height $10^{21}$ (blue) and $10^{22} (red)$}
		\label{fig:Distributions_at_10_22}
\end{figure}

It is interesting to see if this behaviour still continues when going further in the zeta zeros. We use here zeta zeros calculated by Odlyzko at the heights $10^{21}$ and $10^{22}$. Because the series of zeros are quite small, only 10000 zeros (exactly 9799 differences), and consequently statistics not very good, we use Johnson distribution fit for the $delta$-distributions. Figure \ref{fig:Distributions_at_10_22}a shows the $delta$-distributions of zeros at height $10^{22}$. The distributions are still lower at the sites of first zeros (only three zeros seen). Figure \ref{fig:Distributions_at_10_22}b shows the corresponding variances for distributions in \ref{fig:Distributions_at_10_22}a, and in addition variances for $delta$-distributions for the zeros at height $10^{21}$.

\newpage

\section{Conclusions}

We have studied $delta$-distributions of zeros of Riemann zeta function at heights 1, 10 and 100 billionth zero such that we calculate distribution for each difference, $delta(n)$, separately. We used 5 million $deltas$ for these analyses, and showed that statistical properties are very similar for all intervals. The distributions are skewed towards nearest zeta zero, and have local maximum variance and local minimum kurtosis when passing zeta zero. This means that zeta zero $deltas$ are repelling zeros themselves. This property seems to last, at least, until the height $10^{22}$ of zeta zero. 
We, furthermore, have shown that Johnson PDFs SB and SU do fit nicely to the $delta$-distributions, although skewness and kurtosis are changing quite a lot. This indicates the flexibility of Johnson distribution analysis in this kind of statistical application.
\newline
\newline

\textbf{Acknowledgements:}

We acknowledge LMFDB and A.M. Odlyzko for the zeta zero data.


\newpage

\begin{thebibliography}{5}

\bibitem{Cugerone}
Cugerone K. and De Michele C.,
\newblock Johnson SB as general functional form for raindrop size distribution,
\newblock \textit{Water Resources Research}, 51, pp. 6276-6289, 2015.

\bibitem{Johnson}
Johnson N.L.,
\newblock Systems of frequency curves generated by methods of translation,
\newblock \textit{Biometrika}, 36, pp. 149-176, 1949.

\bibitem{Montgomery}
Montgomery H.L.,
\newblock The pair correlation of zeros of the zeta function,
\newblock \textit{Analytic number theory (Proc. Sympos. Pure Math., Vol. XXIV, St. Louis Univ., St.
	Louis, Mo., 1972)}, 1973.

\bibitem{Odlyzko}
Odlyzko, A.M.,
\newblock On the Distribution of Spacings Between Zeros of the Zeta Function,
\newblock \textit{Mathematics of Computation}, 48 (177), pp. 273-308, 1987.

\bibitem{Perez-Marco}
Perez-Marco R.,
\newblock Statistics of Riemann zeros ,
\newblock \textit{arXiv:1112.0346v1 [math.NT]}, 2011.

\bibitem{Snaith}
Snaith N.C.,
\newblock Riemann Zeros and Random Matrix Theory,
\newblock \textit{Milan J. Math.}, 78, pp. 135-152, 2010.

\bibitem{Wheeler}
Wheeler R.E,
\newblock Quantile estimators of Johnson curve parameters,
\newblock \textit{Biometrika}, 67 (3), pp. 725-728, 1980.

\end{thebibliography}
\end{document}